# FOCUSING AT A POINT WITH CAUSTIC CROSSING FOR A CLASS OF NONLINEAR EQUATIONS

RÉMI CARLES AND DAVID LANNES

ABSTRACT. We consider the asymptotic behavior of solutions to nonlinear partial differential equations in the limit of short wavelength. For initial data which cause focusing at one point, we highlight critical indexes as far as the influence of the nonlinearity in the neighborhood of the caustic is concerned. Our results generalize some previous ones proved by the first author in the case of nonlinear Schrödinger equations. We apply them to Hartree, Klein-Gordon and wave equations.

## 1. INTRODUCTION

The mathematical theory of geometrical optics describes the asymptotic behavior of solutions to problems of the form, for instance,

$$\begin{cases} L(\varepsilon,\partial)u^\varepsilon = F(u^\varepsilon), \\ u^\varepsilon_{|t=0} = f(x)e^{i\varphi(x)/\varepsilon}, \end{cases}$$

when the parameter $\varepsilon$ goes to zero. One seeks in general a solution of the form

$$u^\varepsilon(t,x) \underset{\varepsilon \to 0}{\sim} \varepsilon^J U_0\left(t, x, \frac{\phi(t,x)}{\varepsilon}\right),$$

for some profile $U_0$ independent of $\varepsilon$. The phase $\phi$ solves an eikonal equation, and may develop singularities in finite time: a caustic appears, where the asymptotic expansion of geometrical optics ceases to be valid. This phenomenon is well understood for linear equations. The caustic crossing is described by the Maslov index ([27], [12]). For a nonlinear equation, few results are available. Formal computations by Hunter and Keller ([21]) suggest that there exist two distinct notions of critical index, depending on the equations, the nonlinearity, the amplitude of the initial datum and the geometry of propagation. In the above example, this would yield a value $J_p$ such that if $J > J_p$, the nonlinear term is negligible outside the caustic (linear geometric optics is valid), while if $J = J_p$, the solution is described by nonlinear geometric optics outside the caustic. The other value $J_c$ leads to a similar discussion near the caustic (as opposed to outside the caustic). Such distinctions were proved by Joly, Métivier and Rauch ([24]) for some nonlinear wave equations, and by the first author ([3]) for some nonlinear Schrödinger equations.

The aim of this paper is to provide some generalization of the main result of [3]. In this article, it was proved that for some semi-classical nonlinear Schrödinger equations with quadratically oscillating initial data, the solution focuses at one

---

point, and the caustic crossing is described at leading order by a scattering operator. Consider the initial value problem,

$$(1.1) \quad \begin{cases} i\varepsilon\partial_t u^\varepsilon + \dfrac{1}{2}\varepsilon^2 \Delta u^\varepsilon = \varepsilon^{n\sigma} |u^\varepsilon|^{2\sigma} u^\varepsilon, & (t,x) \in \mathbb{R}_+ \times \mathbb{R}^n, \\ u^\varepsilon_{|t=0} = e^{-i\frac{|x|^2}{2\varepsilon}} f(x), \end{cases}$$

with $\varepsilon \in ]0,1]$, and suitable assumptions on $\sigma$. We assume that $f \in \Sigma := H^1 \cap \mathcal{F}(H^1)$. The asymptotic behavior of the solution $u^\varepsilon$ is studied, as $\varepsilon$ goes to zero: quadratic oscillations cause focusing at the origin at time $t = 1$, and the scaling of the nonlinearity (the presence of the factor $\varepsilon^{n\sigma}$) makes the influence of the right-hand side of (1.1) negligible away from the focal point ("linear propagation", see [3]). On the other hand, the caustic crossing takes some nonlinear effects into account ("nonlinear caustic"), and is described at leading order by the (nonlinear) scattering operator associated with

$$(1.2) \quad i\partial_t \psi + \frac{1}{2}\Delta \psi = |\psi|^{2\sigma}\psi.$$

Explicitly, the following asymptotics holds in $\Sigma$,

$$(1.3) \quad u^\varepsilon(t,x) = \begin{cases} \dfrac{1}{(1-t)^{n/2}} f\left(\dfrac{x}{1-t}\right) e^{i\frac{|x|^2}{2\varepsilon(t-1)}} + o(1), & \text{if } t < 1, \\ \dfrac{e^{-in\frac{\pi}{2}}}{(t-1)^{n/2}} Zf\left(\dfrac{x}{1-t}\right) e^{i\frac{|x|^2}{2\varepsilon(t-1)}} + o(1), & \text{if } t > 1, \end{cases}$$

where

$$Z = \mathcal{F}^{-1} \circ S \circ \mathcal{F},$$

$S$ is the nonlinear scattering operator associated to (1.2) (see e.g. [11]), and

$$\widehat{f}(\xi) = \mathcal{F}f(\xi) := \frac{1}{(2\pi)^{n/2}} \int_{\mathbb{R}^n} e^{-ix\cdot\xi} f(x) dx.$$

Two phenomena are to be noticed is the above description. First, the phase shift $e^{-in\frac{\pi}{2}}$ after the focus crossing is the usual Maslov index ([12]). It can be understood as a linear phenomenon. The other point is the change of profile. Before focusing, the profile $f$ describes the behavior of $u^\varepsilon$. Past the caustic, we have another profile, $Zf$, which is deduced from $f$ by a nonlinear manipulation (essentially a nonlinear scattering operator). We prove in this paper that these two phenomena occur for a class of equations, such as Hartree, Klein-Gordon and the wave equations.

The above asymptotics actually stems from a stronger approximation, as shown in [3]. Let $W_-$ denote the wave operator in $-\infty$ associated to (1.2). Roughly speaking, for $\psi_-$ an asymptotic state, $W_-\psi_- = \psi_{|t=0}$ where $\psi$ solves (1.2) with $\psi_-$ as initial data prescribed at $-\infty$,

$$e^{-i\frac{t}{2}\Delta}\psi(t,x)\big|_{t=-\infty} = \psi_-(x).$$

If $\psi_-$ is defined as $\psi_- = \mathcal{F}f$, and $\psi$ is given by $W_-\psi_- = \psi_{|t=0}$, then the following asymptotics holds in $L^\infty_{loc}(\mathbb{R}^+_t, \Sigma)$,

$$(1.4) \quad u^\varepsilon(t,x) \underset{\varepsilon \to 0}{\sim} \frac{1}{\varepsilon^{n/2}} \psi\left(\frac{t-1}{\varepsilon}, \frac{x}{\varepsilon}\right).$$

In particular, at the focusing time,

$$(1.5) \qquad u^\varepsilon(1,x) \underset{\varepsilon \to 0}{\sim} \frac{1}{\varepsilon^{n/2}} \varphi\left(\frac{x}{\varepsilon}\right),$$

where $\varphi = \psi_{|t=0}$. Conversely, assume that (1.5) holds for some function $\varphi$, and that $u^\varepsilon$ solves

$$i\varepsilon \partial_t u^\varepsilon + \frac{1}{2}\varepsilon^2 \Delta u^\varepsilon = \varepsilon^{n\sigma} |u^\varepsilon|^{2\sigma} u^\varepsilon.$$

Splitting the variables around the point $(t,x) = (1,0)$ at scale $\varepsilon$,

$$u^\varepsilon(t,x) = \frac{1}{\varepsilon^{n/2}} \psi^\varepsilon\left(\frac{t-1}{\varepsilon}, \frac{x}{\varepsilon}\right),$$

we see that $\psi^\varepsilon$ solves (1.2). The continuity of the wave operator $W_-$ and the global well-posedness of (1.2) insure that if $\sigma$ is sufficiently large, and $\varphi$ sufficienlty smooth and/or small, then (1.4) holds in $L^\infty_{loc}(\mathbb{R}^+_t, \Sigma)$, where $\psi_{|t=0} = \varphi$.

When global well-posedness is known for (1.2), then it is also known that the wave operators $W_\pm$ are invertible, and

$$\psi(T,X) \underset{T \to \pm\infty}{\sim} e^{i\frac{T}{2}\Delta} \psi_\pm(X),$$

where $\psi_\pm = W_\pm^{-1} \varphi$. A stationary phase argument (see e.g. [33]) shows that

$$e^{i\frac{T}{2}\Delta} \psi_\pm(X) \underset{T \to \pm\infty}{\sim} \frac{e^{\mp in\pi/4}}{|T|^{n/2}} e^{i\frac{X^2}{2T}} \mathcal{F}^{-1} \psi_\pm\left(-\frac{X}{T}\right),$$

and the substitution

$$(T,X) = \left(\frac{t-1}{\varepsilon}, \frac{x}{\varepsilon}\right)$$

yields the asymptotics (1.3), with $f = e^{in\pi/4} \mathcal{F}^{-1} \psi_-$. Thus, if we start from (1.5), then proving (1.3) can be seen as a byproduct of the scattering theory for (1.2) and of the asymptotics of the free evolution group $e^{it\Delta}$. This yields a simple proof of the case "linear propagation, nonlinear caustic" in [3]. Similarly, the case 'linear propagation, linear caustic" can be regarded as a consequence of the continuity of the scattering operator $S$ at the origin. That case corresponds to the situation

$$(1.6) \qquad \begin{cases} i\varepsilon \partial_t u^\varepsilon + \dfrac{1}{2}\varepsilon^2 \Delta u^\varepsilon = \varepsilon^\alpha |u^\varepsilon|^{2\sigma} u^\varepsilon, & (t,x) \in \mathbb{R}_+ \times \mathbb{R}^n, \\ u^\varepsilon_{|t=0} = e^{-i\frac{|x|^2}{2\varepsilon}} f(x), \end{cases}$$

with $\alpha > n\sigma$. The conclusion stated in [3] is that $u^\varepsilon$ behaves asymptotically like the solution of the linear equation

$$(1.7) \qquad \begin{cases} i\varepsilon \partial_t v^\varepsilon + \dfrac{1}{2}\varepsilon^2 \Delta v^\varepsilon = 0, & (t,x) \in \mathbb{R}_+ \times \mathbb{R}^n, \\ v^\varepsilon_{|t=0} = e^{-i\frac{|x|^2}{2\varepsilon}} f(x). \end{cases}$$

This means that the nonlinearity in (1.6) is negligible at leading order, away from the focus, as well as near the focus. For a new proof of this result, introduce the scaling

$$u^\varepsilon(t,x) = \varepsilon^{-\alpha/2\sigma} \psi^\varepsilon\left(\frac{t-1}{\varepsilon}, \frac{x}{\varepsilon}\right).$$

Then $\psi^\varepsilon$ solves (1.2). Denote $U_0(t) = e^{i\frac{t}{2}\Delta}$. Then one can compute, in $\Sigma$,

$$U_0\left(\frac{1}{\varepsilon}\right)\psi^\varepsilon\left(\frac{-1}{\varepsilon}\right) = \varepsilon^{\alpha/2\sigma - n/2}\left(\widehat{f} + o(1)\right), \quad \text{as } \varepsilon \to 0.$$

Since the power of $\varepsilon$ in the right hand side is positive, this means that we consider a *small* asymptotic state. Then the global well-posedness of (1.2) for small data proves that

$$U_0(-t)\psi^\varepsilon(t,x) \underset{\varepsilon \to 0}{\sim} \varepsilon^{\alpha/2\sigma - n/2}\psi_-, \quad \text{in } L^\infty(\mathbb{R}, \Sigma),$$

which yields the asymptotics

$$\|U_0(-t)(u^\varepsilon - v^\varepsilon)(t, \cdot)\|_{L^\infty(\mathbb{R},\Sigma)} = o(1).$$

In this paper, we repeat this argument to a class of nonlinear equations, for which well-posedness, existence of wave operators and their asymptotic completeness, along with asymptotics for the associated linear equation are known. This includes in particular the case of Hartree, Klein-Gordon and semilinear wave equations. In Section 2, we describe this general approach and state our main results (Theorems 2.5 and 2.7). Our first two applications concern Hartree and Klein-Gordon equations in Sections 3 and 4 respectively. In Section 5, we analyze the semilinear wave equation, and compare our results with some proved in [1], [13], [10] and [9]. Finally, we conclude this paper with some remarks and comments in Section 6.

## 2. The general approach

In this section, we describe an algorithm which generalizes the approach presented in introduction.

2.1. **A linear equation.** We consider here a linear system of $q$ partial differential equations

$$(2.1) \quad \begin{cases} L(\partial)\psi = 0, \\ \psi_{|t=0} = \varphi. \end{cases}$$

The operator $L(\partial)$ is a differential operator in the variables $(t,x) \in \mathbb{R} \times \mathbb{R}^n$, with constant coefficients, acting on functions $\psi : (t,x) \in \mathbb{R} \times \mathbb{R}^n \to \mathbb{C}^q$. We assume that it is of the form

$$L(\partial) = \partial_t + p(\partial_{x_1}, \ldots, \partial_{x_n}),$$

where $p(\cdot)$ is a polynomial whose coefficient are matrices in $\mathcal{M}_q(\mathbb{C})$.

*Remark* 2.1. The reason why we assume that $L(\partial)$ has constant coefficients is geometric. Non constant coefficients may lead to the phenomenon of refocusing (see for instance [22]), which we want to exclude.

We suppose that $\varphi \in X$, where $X$ is a Banach space invariant with respect to dilations; if $\varphi \in X$, then for any $\lambda > 0$, $\varphi(\lambda \cdot) \in X$.

**Assumption 1.** *We assume that the following properties hold:*
- *The Cauchy problem (2.1) has a unique solution $\psi \in C(\mathbb{R}; X)$, and there exists a group $U_0$, unitary on $X$, such that $\psi(t) = U_0(t)\varphi$, for all $t \in \mathbb{R}$.*

- The solutions are asymptotically of size $|t|^{-\alpha}$. There exists $\alpha \geq 0$, such that for every $\varphi \in \mathcal{S}(\mathbb{R}^n; \mathbb{C}^q)$, there exist $c_\pm = c_\pm(\varphi) < \infty$ such that
$$\limsup_{t \to \pm\infty} |t|^\alpha \|U_0(t)\varphi(\cdot)\|_{L^\infty(\mathbb{R}^n, \mathbb{C}^q)} = c_\pm,$$
and there exist $\varphi_\pm$ such that $c_\pm$ are nonzero.
- This group exhibits an asymptotic behavior in $\pm\infty$. There exist mappings $A_\pm : X \to C(\mathbb{R}; X)$ such that for every $\varphi \in X$,

(2.2) $$\left\| \varphi - U_0(-t) \left( \frac{1}{|t|^\alpha} A_\pm(\varphi)(t, \cdot) \right) \right\|_X \xrightarrow[t \to \pm\infty]{} 0 .$$

*Remark* 2.2. The existence of finite constants $c_\pm$ which are not identically equal to zero ensures the uniqueness of $\alpha$ in the above definition. Moreover, it means that the initial data we consider in the initial value problem (2.5) below are *morally* of order $O(1)$ in $L^\infty$ (in general, nothing allows us to say that they are in $L^\infty$).

*Example.* In the case of the Schrödinger equation, we can take $X = L^2(\mathbb{R}^n)$. We have $U_0(t) = e^{i\frac{t}{2}\Delta}$, which is unitary on $H^s(\mathbb{R}^n)$ for every $s$, and for every $\varphi \in L^2$ (see e.g. [33]),
$$\left\| U_0(t)\varphi - e^{-in\frac{\pi}{4}\operatorname{sgn} t} \frac{e^{i\frac{x^2}{2t}}}{t^{n/2}} \widehat{\varphi}\left(\frac{x}{t}\right) \right\|_{L^2} \xrightarrow[t \to \pm\infty]{} 0.$$
Therefore, the above assumptions are satisfied, with $\alpha = n/2$, and
$$A_\pm(\varphi)(t, x) = e^{\mp in\frac{\pi}{4}} e^{i\frac{x^2}{2t}} \widehat{\varphi}\left(\frac{x}{t}\right).$$

*Remark* 2.3. For nonlinear Schrödinger equations, the above assumption is not exactly verified in the case analyzed in [3]. Indeed, the group $U_0(t) = e^{i\frac{t}{2}\Delta}$ is unitary on $H^s(\mathbb{R}^n)$ for any $s$, but not on $\mathcal{F}(H^1)$, hence not on $\Sigma$, which is the space $X$ considered in [3]. The assumption that $U_0$ is unitary on $X$ eases the readability of the statements and proofs below, and is satisfied in most of the other examples we consider, this is the reason why we make it.

### 2.2. A nonlinear equation.
We now consider the semilinear equation

(2.3) $$\begin{cases} L(\partial)\psi = F(x, \psi), \\ \psi_{|t=0} = \varphi. \end{cases}$$

We assume that $F$ is a nonlinearity which is $p$-homogeneous in its last variable for some $p > 1$,
$$F(x, \lambda z) = \lambda^p F(x, z), \quad \forall (x, \lambda, z) \in \mathbb{R}^n \times \mathbb{R}_+^* \times \mathbb{C}^q.$$

*Remark* 2.4. We assume that the nonlinearity depends on the space variable $x$ to compare with the results of Carles and Rauch (see Section 5), and to apply our result to the Hartree equation (Section 3).

We suppose there exists a (short range) scattering theory in the following sense. For $t_0 \in \mathbb{R}_-$ and $\psi_- \in X$, consider the initial value problem

(2.4) $$\begin{cases} L(\partial)\psi = F(x, \psi), \\ U_0(-t)\psi(t)_{|t=t_0} = \psi_-. \end{cases}$$

**Assumption 2.** *For every $\psi_- \in X$, there exists $T < \infty$ such that*

- For every $t_0 \in [-\infty, -T]$, (2.4) has a unique solution in $C(]-\infty, -T]; X)$.
- The solution $\psi \in C(]-\infty, -T]; X)$ depends continuously on $\psi_- \in X$ and $t_0 \in [-\infty, -T]$.
- Define $\widetilde{\psi}(t) := U_0(-t)\psi(t)$. Then $\widetilde{\psi} \in C([-\infty, -T]; X)$. In particular, if $t_0 = -\infty$, $\lim\limits_{t \to -\infty} \widetilde{\psi}(t) = \psi_-$.

When $\psi$ is defined up to time $t = 0$, the map $W_- : \psi_- \mapsto \psi_{|t=0}$ is the wave operator at $-\infty$. The continuity assumption (second point of Assumption 2) is in general a by-product of the existence of the wave operator $W_-$, obtained by a fixed point argument.

We also need the asymptotic completeness of the wave operators, as well as global well-posedness.

**Assumption 3.** *We suppose that the following properties hold.*

- *The initial value problem* (2.3) *is globally well-posed in* $C(\mathbb{R}; X)$; *for every $R > 0$ and every $\varphi \in B_R := \{\varphi \in X \ ; \ \|\varphi\|_X < R\}$, the Cauchy problem* (2.3) *has a unique solution $\psi \in C(\mathbb{R}; X)$. Moreover, the map $\varphi \mapsto \widetilde{\psi}$ is uniformly continuous from $B_R$ to $C(\mathbb{R}; X)$.*
- *Let $\varphi \in X$, and $\psi$ the solution of* (2.3)*. There exist $\psi_\pm \in X$ such that*

$$\left\|\widetilde{\psi}(t) - \psi_\pm\right\|_X = \|U_0(-t)\psi(t) - \psi_\pm\|_X \xrightarrow[t \to \pm\infty]{} 0.$$

**Definition.** When Assumptions 2 and 3 and satisfied, we define the wave operators at $\pm\infty$ as $W_\pm : \psi_\pm \mapsto \psi_{|t=0}$, and the scattering operator $S : \psi_- \mapsto \psi_+$.

2.3. **A general class of focusing initial value problems.** In the case of the Schrödinger equation, the scattering operator associated to (1.2) is used to describe the behavior of the initial value problem (1.1). In the general framework considered in this paper, the properties of the generic nonlinear equation (2.3) are used to describe the solution $u^\varepsilon$ of the initial value problem

$$(2.5) \quad \begin{cases} L(\varepsilon\partial)u^\varepsilon(t,x) = \varepsilon^{\alpha(p-1)} F\left(\dfrac{x}{\varepsilon}, u^\varepsilon(t,x)\right), \\ u^\varepsilon(0,x) = A_-(\psi_-)\left(\dfrac{-1}{\varepsilon}, \dfrac{x}{\varepsilon}\right) + r^\varepsilon(x), \end{cases}$$

where $r^\varepsilon \in X$ is asymptotically small in a sense we make precise below. Such initial data cause focusing in the linear case. As $\varepsilon \to 0$, the solution to the problem

$$L(\varepsilon\partial)w^\varepsilon = 0 \ , \quad w^\varepsilon(0,x) = A_-(\psi_-)\left(\frac{-1}{\varepsilon}, \frac{x}{\varepsilon}\right),$$

is given by

$$w^\varepsilon(t,x) = U_0\left(\frac{t}{\varepsilon}\right) A_-(\psi_-)\left(\frac{-1}{\varepsilon}, \frac{x}{\varepsilon}\right).$$

From Assumption 1, we have

$$w^\varepsilon(t,x) \underset{\varepsilon \to 0}{\sim} \frac{1}{\varepsilon^\alpha} U_0\left(\frac{t-1}{\varepsilon}\right) \psi_-\left(\frac{x}{\varepsilon}\right).$$

Therefore, when $\alpha > 0$, the amplitude of $v^\varepsilon$ grows at time $t = 1$, from $O(1)$ initially, to $O(\varepsilon^{-\alpha})$; it focuses at the origin at time $t = 1$.

*Example.* In the case of the nonlinear Schrödinger equation, $\psi$ solves (1.2). The nonlinearity is $(2\sigma+1)$-homogeneous, and we saw that $\alpha = n/2$, so $u^\varepsilon$ solves
$$i\varepsilon \partial_t u^\varepsilon + \frac{1}{2}\varepsilon^2 \Delta u^\varepsilon = \varepsilon^{n\sigma}|u^\varepsilon|^{2\sigma}u^\varepsilon.$$
We also noticed that $A_-(\psi_-)(t,x) = e^{i\frac{x^2}{2t}}\widehat{\varphi}\left(\frac{x}{t}\right)$, so
$$u^\varepsilon(0,x) = A_-(\psi_-)\left(\frac{-1}{\varepsilon}, \frac{x}{\varepsilon}\right) = e^{-i\frac{x^2}{2\varepsilon}}\widehat{\varphi}(-x).$$
We thus retrieve the problem presented in introduction, with $f(x) := \widehat{\varphi}(-x)$.

In general, the asymptotic behavior of $u^\varepsilon$ is not exactly given with respect to the norm of $X$. For $\varphi \in X$ and $\varepsilon \in ]0,1]$, define
$$\|\varphi\|_{X,\varepsilon} := \|\varepsilon^\alpha \varphi(\varepsilon\cdot)\|_X. \tag{2.6}$$
This norm $\|\cdot\|_{X,\varepsilon}$ is the appropriate one to describe the behavior of $u^\varepsilon$.

*Example.* In the case the nonlinear Schrödinger equation, for $X = H^1(\mathbb{R}^n)$, this definition yields
$$\|\varphi\|_{H^1,\varepsilon} = \|\varphi\|_{L^2} + \|\varepsilon \nabla \varphi\|_{L^2},$$
and the dependence of this norm upon $\varepsilon$ appears in the presence of $\varepsilon$-derivatives, instead of the classical derivatives. One must not be surprised: when oscillations with a short wavelength $\varepsilon$ are considered, such norm have proved to be efficient in geometric optics. They were introduced by Guès in [17] (see also [23], [34]) to justify high-frequency approximations in nonlinear problems.

To make the meaning of (2.5) definitely clear, we state it in a slightly different fashion, that will remove any ambiguity when it turns to study the wave equation (or more generally, an equation of order at least two written as a vector-valued first order equation). Consider the initial value problem
$$\begin{cases} L(\partial)\psi^\varepsilon(t,x) = F(x, \psi^\varepsilon(t,x)), \\ \psi^\varepsilon\left(-\frac{1}{\varepsilon}, x\right) = \varepsilon^\alpha A_-(\psi_-)\left(\frac{-1}{\varepsilon}, x\right) + \varepsilon^\alpha r^\varepsilon(\varepsilon x). \end{cases} \tag{2.7}$$
Notice that the function $\psi^\varepsilon$ depends on $\varepsilon$ because the time at which its value is prescribed, as well as this initial datum, depend on $\varepsilon$. Now define $u^\varepsilon$ as
$$u^\varepsilon(t,x) := \frac{1}{\varepsilon^\alpha}\psi^\varepsilon\left(\frac{t-1}{\varepsilon}, \frac{x}{\varepsilon}\right). \tag{2.8}$$
Since the nonlinearity $F$ is $p$-homogeneous in its last argument, it is easy to check that (2.7) is equivalent to (2.5), through (2.8).

2.4. **The main results.** In this paragraph, we prove that the caustic crossing for $u^\varepsilon$ as given by (2.5) leads to the same phenomenon as in [3]. We first treat the case where the focusing is described by the nonlinear scattering operator $S$, and show next that the power $\alpha(p-1)$ in (2.5) is critical: if it is replaced by $\delta > \alpha(p-1)$, then the nonlinearity is globally negligible at leading order.

We first state our main result.

**Theorem 2.5.** *Let $\psi_-, r^\varepsilon \in X$, and suppose that Assumptions 1, 2 and 3 are satisfied. If $r^\varepsilon$ is asymptotically small, $\|r^\varepsilon\|_{X,\varepsilon} \xrightarrow[\varepsilon \to 0]{} 0$, then the following holds.*

- There exists $\varepsilon_0 > 0$ such that for all $0 < \varepsilon < \varepsilon_0$, there exists a unique solution $u^\varepsilon \in C(\mathbb{R}_+, X)$ to the Cauchy problem (2.5).
- For every $0 < \tau < 1$, $u^\varepsilon$ satisfies the following asymptotics,

(2.9) $$\sup_{t \in [0,\tau]} \left\| u^\varepsilon(t, \cdot) - \frac{1}{(1-t)^\alpha} A_-(\psi_-)\left(\frac{t-1}{\varepsilon}, \frac{\cdot}{\varepsilon}\right) \right\|_{X,\varepsilon} \xrightarrow[\varepsilon \to 0]{} 0 \ .$$

- For every $\tau > 1$, $u^\varepsilon$ satisfies the following asymptotics,

(2.10) $$\sup_{t \geq \tau} \left\| u^\varepsilon(t, \cdot) - \frac{1}{(t-1)^\alpha} A_+(\psi_+)\left(\frac{t-1}{\varepsilon}, \frac{\cdot}{\varepsilon}\right) \right\|_{X,\varepsilon} \xrightarrow[\varepsilon \to 0]{} 0 \ ,$$

where $\psi_+ = S\psi_-$.
- There exists a "caustic profile". Let $\varphi = W_-\psi_- \in X$. Then

$$\left\| u^\varepsilon(1, \cdot) - \frac{1}{\varepsilon^\alpha} \varphi\left(\frac{\cdot}{\varepsilon}\right) \right\|_{X,\varepsilon} \xrightarrow[\varepsilon \to 0]{} 0.$$

*Proof.* **Step 1. Global existence of $u^\varepsilon$.** Define $\psi^\varepsilon$ by (2.8), that is

$$u^\varepsilon(t,x) = \frac{1}{\varepsilon^\alpha} \psi^\varepsilon\left(\frac{t-1}{\varepsilon}, \frac{x}{\varepsilon}\right).$$

Then $u^\varepsilon$ solves (2.5) if and only if $\psi^\varepsilon$ solves (2.7), or, equivalently

(2.11) $$\begin{cases} L(\partial)\psi^\varepsilon = F(x, \psi^\varepsilon) \ , \\ U_0\left(\frac{1}{\varepsilon}\right)\psi^\varepsilon\left(-\frac{1}{\varepsilon}\right) = \varepsilon^\alpha U_0\left(\frac{1}{\varepsilon}\right) A_-(\psi_-)\left(-\frac{1}{\varepsilon}\right) + \varepsilon^\alpha U_0\left(\frac{1}{\varepsilon}\right) r^\varepsilon(\varepsilon \cdot) \ . \end{cases}$$

From Assumption 2, we know that there exists $T$ finite such that for all $t_0 \in [-\infty, -T]$ and $\psi_{-,t_0} \in X$, the Cauchy problem (2.4) with datum $U_0(-t)\psi(t)|_{t=t_0} = \psi_{-,t_0}$ admits a unique solution in $C(]-\infty, -T]; X)$, and that this solution depends continuously on $t_0$ and $\psi_{-,t_0}$. Since by Assumption 1,

$$\left\| \varepsilon^\alpha U_0\left(\frac{1}{\varepsilon}\right) A_-(\psi_-)\left(-\frac{1}{\varepsilon}\right) - \psi_- \right\|_X \xrightarrow[\varepsilon \to 0]{} 0 \ ,$$

and since we supposed that $\|\varepsilon^\alpha r^\varepsilon(\varepsilon \cdot)\|_X \to 0$ as $\varepsilon \to 0$, this yields that the solution $\psi^\varepsilon$ to (2.11) is well-defined in $C(]-\infty, T]; X)$ for some $T > 0$ *independent of $\varepsilon$* (provided that $\varepsilon$ is small enough). Moreover, one has the following estimate,

(2.12) $$\sup_{t \leq -T} \left\| \psi^\varepsilon\left(\frac{t-1}{\varepsilon}, \cdot\right) - \psi^*\left(\frac{t-1}{\varepsilon}, \cdot\right) \right\|_X \xrightarrow[\varepsilon \to 0]{} 0 \ ,$$

where $\psi^*$ is defined as

$$\begin{cases} L(\partial)\psi^* = F(x, \psi^*) \ , \\ U_0(-t)\psi^*(t)|_{t=-\infty} = \psi_- \ . \end{cases}$$

The fact that $\psi^\varepsilon$ is globally defined in $C(\mathbb{R}; X)$ is then a consequence of Assumption 3. As already said, existence and uniqueness of $u^\varepsilon$ follow from existence and uniqueness of $\psi^\varepsilon$.

**Step 2. Estimates before the caustic.** Let $0 < \tau < 1$ and $\psi$ be the solution of (2.3) with datum $\varphi := W_-\psi_-$. Define $v^\varepsilon$ by

$$v^\varepsilon(t,x) = \frac{1}{\varepsilon^\alpha} \psi\left(\frac{t-1}{\varepsilon}, \frac{x}{\varepsilon}\right).$$

Then one has

(2.13) $$\sup_{t\in[0,\tau]} \|u^\varepsilon(t,\cdot) - v^\varepsilon(t,\cdot)\|_{X,\varepsilon} \to 0, \quad \text{as} \quad \varepsilon \to 0.$$

Indeed, it follows from the definition of $W_-$ that $\psi = \psi^*$ where $\psi^*$ is defined in Step 1. Therefore

$$\begin{aligned} \|u^\varepsilon(t,\cdot) - v^\varepsilon(t,\cdot)\|_{X,\varepsilon} &= \left\|\psi\left(\frac{t-1}{\varepsilon},\cdot\right) - \psi^\varepsilon\left(\frac{t-1}{\varepsilon},\cdot\right)\right\|_X \\ &= \left\|\psi^*\left(\frac{t-1}{\varepsilon},\cdot\right) - \psi^\varepsilon\left(\frac{t-1}{\varepsilon},\cdot\right)\right\|_X. \end{aligned}$$

Since for $t \in [0,\tau]$ and $\varepsilon$ small enough, $\frac{t-1}{\varepsilon} \leq \frac{\tau-1}{\varepsilon} \leq -T$, (2.13) is a consequence of (2.12).

**Step 3. Asymptotics before the caustic.** We have

(2.14) $$\sup_{t\in[0,\tau]} \left\|u^\varepsilon(t,\cdot) - \frac{1}{(1-t)^\alpha} A_-(\psi_-)\left(\frac{t-1}{\varepsilon}, \frac{\cdot}{\varepsilon}\right)\right\|_{X,\varepsilon} \to 0 \quad \text{as} \quad \varepsilon \to 0.$$

Indeed, we have

$$\begin{aligned} & \left\|v^\varepsilon(t,\cdot) - \frac{1}{(1-t)^\alpha} A_-(\psi_-)\left(\frac{t-1}{\varepsilon}, \frac{\cdot}{\varepsilon}\right)\right\|_{X,\varepsilon} \\ &= \left\|\psi\left(\frac{t-1}{\varepsilon}\right) - \frac{\varepsilon^\alpha}{(1-t)^\alpha} A_-(\psi_-)\left(\frac{t-1}{\varepsilon},\cdot\right)\right\|_X \\ &= \left\|U_0\left(-\frac{t-1}{\varepsilon}\right)\psi\left(\frac{t-1}{\varepsilon}\right) - U_0\left(-\frac{t-1}{\varepsilon}\right)\frac{\varepsilon^\alpha}{(1-t)^\alpha} A_-(\psi_-)\left(\frac{t-1}{\varepsilon},\cdot\right)\right\|_X \\ &\leq \left\|U_0\left(-\frac{t-1}{\varepsilon}\right)\psi\left(\frac{t-1}{\varepsilon}\right) - \psi_-\right\|_X \\ &\quad + \left\|U_0\left(-\frac{t-1}{\varepsilon}\right)\frac{\varepsilon^\alpha}{(1-t)^\alpha} A_-(\psi_-)\left(\frac{t-1}{\varepsilon},\cdot\right) - \psi_-\right\|_X. \end{aligned}$$

From Assumptions 1 and 3, the two terms of the r.h.s. of the above inequality go to 0 as $\varepsilon$ goes to 0, uniformly for $t \in [0,\tau]$. Step 2 and a triangle inequality then yield (2.14).

**Step 4. At the caustic and beyond.** The methods are exactly the same as before the caustic. Global-posedness for (2.3) implies that we can take the supremum in time in (2.13) on $\mathbb{R}_+$, that is,

(2.15) $$\sup_{t\geq 0} \|u^\varepsilon(t,\cdot) - v^\varepsilon(t,\cdot)\|_{X,\varepsilon} \to 0, \quad \text{as} \quad \varepsilon \to 0.$$

In particular, taking $t = 1$ yields the last point of the theorem. Moreover, (2.14) remains true for $t \geq \tau > 1$, provided that $A_-(\psi_-)$ is replaced by $A_+(\psi_+)$. $\square$

*Remark* 2.6. In the above theorem, Assumption 3 is not needed to prove the second point. In other words, existence and asymptotics for $0 \leq t \leq \tau < 1$ still hold if we suppose only that Assumptions 1 and 2 are satisfied.

We now prove that the scaling of the nonlinearity in (2.5) is critical in terms of the power of $\varepsilon$. Let $\delta > \alpha(p-1)$, and suppose that $u^\varepsilon$ solves

(2.16)
$$\begin{cases} L(\varepsilon\partial)u^\varepsilon = \varepsilon^\delta F\left(\dfrac{x}{\varepsilon}, u^\varepsilon\right), \\ u^\varepsilon(0,x) = A_-(\psi_-)\left(\dfrac{-1}{\varepsilon}, \dfrac{x}{\varepsilon}\right) + r^\varepsilon(x). \end{cases}$$

We prove in the next theorem that the nonlinearity in (2.16) is not relevant in the sense that $u^\varepsilon$ behaves asymptotically as the solution $w^\varepsilon$ of the free equation

(2.17)
$$\begin{cases} L(\varepsilon\partial)w^\varepsilon = 0, \\ w^\varepsilon(0,x) = A_-(\psi_-)\left(\dfrac{-1}{\varepsilon}, \dfrac{x}{\varepsilon}\right). \end{cases}$$

Such a result requires an additional assumption on the small data case, whose relevance is discussed in Section 2.6.

**Assumption 4.** *There exist $T < \infty$ such that for any $t_0 \in [-\infty, -T]$ and any $\psi_-$ small enough in $X$, the Cauchy problem (2.4) has a unique solution $\psi \in C(\mathbb{R}; X)$. Moreover, there exists $C = C(\|\psi_-\|_X)$ independent of $t_0 \in [-\infty, -T]$ such that*

$$\sup_{t > t_0} \|U_0(-t)\psi(t,\cdot) - \psi_-\|_X \leq C\left(\|\psi_-\|_X\right)\|\psi_-\|_X, \quad \text{with } C(\lambda) \xrightarrow[\lambda \to 0]{} 0.$$

**Theorem 2.7.** *Let $\psi_-, r^\varepsilon \in X$, and suppose that Assumptions 1 and 4 are satisfied. Let $u^\varepsilon$ be the solution of the Cauchy problem (2.16) with $\delta > \alpha(p-1)$, and $w^\varepsilon$ the solution of (2.17). Assume that $\|r^\varepsilon\|_{X,\varepsilon} \to 0$ as $\varepsilon \to 0$. Then for $\varepsilon$ sufficiently small, $u^\varepsilon \in C(\mathbb{R}_+; X)$ and*

$$\sup_{t \geq 0}\|u^\varepsilon(t,\cdot) - w^\varepsilon(t,\cdot)\|_{X,\varepsilon} \xrightarrow[\varepsilon \to 0]{} 0.$$

*Proof.* Define $\psi^\varepsilon$ by

$$u^\varepsilon(t,x) = \frac{1}{\varepsilon^\beta}\psi^\varepsilon\left(\frac{t-1}{\varepsilon}, \frac{x}{\varepsilon}\right),$$

and choose $\beta$ so that $\psi^\varepsilon$ solves $L(\partial)\psi^\varepsilon = F(x, \psi^\varepsilon)$. This yields $\beta = \dfrac{\delta}{p-1} > \alpha$, and $\psi^\varepsilon$ solves the Cauchy problem

$$\begin{cases} L(\partial)\psi^\varepsilon = F(x, \psi^\varepsilon), \\ U_0(-t)\psi^\varepsilon(t)\big|_{t=-1/\varepsilon} = \varepsilon^{\beta-\alpha}\varepsilon^\alpha U_0\left(\dfrac{1}{\varepsilon}\right)\left(A_-(\psi_-)\left(\dfrac{-1}{\varepsilon}, x\right) + r^\varepsilon(\varepsilon\cdot)\right). \end{cases}$$

We now rescale the solution $w^\varepsilon$ of (2.17) as

$$w^\varepsilon(t,x) = \frac{1}{\varepsilon^\beta}\widetilde{\psi}^\varepsilon\left(\frac{t-1}{\varepsilon}, \frac{x}{\varepsilon}\right).$$

Since $w^\varepsilon$ satisfies a free equation, one has

$$\widetilde{\psi}^\varepsilon\left(\frac{t-1}{\varepsilon}, \cdot\right) = U_0\left(\frac{t}{\varepsilon}\right)\varepsilon^\beta A_-(\psi_-)\left(\frac{-1}{\varepsilon}, \cdot\right).$$

Therefore the difference $u^\varepsilon - w^\varepsilon$ can be estimated as

$$\|u^\varepsilon(t,\cdot) - w^\varepsilon(t,\cdot)\|_{X,\varepsilon} = \left\|\varepsilon^{\alpha-\beta}\psi^\varepsilon\left(\frac{t-1}{\varepsilon}\right) - \varepsilon^{\alpha-\beta}\widetilde{\psi}^\varepsilon\left(\frac{t-1}{\varepsilon}\right)\right\|_X$$

$$= \left\|\varepsilon^{\alpha-\beta}\psi^\varepsilon\left(\frac{t-1}{\varepsilon}\right) - \varepsilon^\alpha U_0\left(\frac{t}{\varepsilon}\right) A_-(\psi_-)\left(\frac{-1}{\varepsilon}\right)\right\|_X$$

Since $U_0$ is unitary on $X$, this quantity is also equal to

$$\left\|\varepsilon^{\alpha-\beta}U_0\left(\frac{1-t}{\varepsilon}\right)\psi^\varepsilon\left(\frac{t-1}{\varepsilon}\right) - \varepsilon^\alpha U_0\left(\frac{1}{\varepsilon}\right) A_-(\psi_-)\left(\frac{-1}{\varepsilon}\right)\right\|_X.$$

From the triangle inequality, it is estimated by

$$\left\|\varepsilon^{\alpha-\beta}U_0\left(\frac{1-t}{\varepsilon}\right)\psi^\varepsilon\left(\frac{t-1}{\varepsilon}\right) - \psi_-\right\| + \left\|\varepsilon^\alpha U_0\left(\frac{1}{\varepsilon}\right) A_-(\psi_-)\left(\frac{-1}{\varepsilon}\right) - \psi_-\right\|_X.$$

The first of the above two terms goes to 0 with $\varepsilon$ as a consequence of Assumption 4, and so does the second one thanks to Assumption 1. □

*Remark* 2.8. The criticality proved above can be understood as the criticality of the size of initial data. With the change of unknown $\widetilde{u}^\varepsilon = \varepsilon^{-\gamma} u^\varepsilon$, $\gamma \geq \alpha$, (2.5) and (2.16) turn into

$$\begin{cases} L(\varepsilon\partial)\widetilde{u}^\varepsilon = F\left(\frac{x}{\varepsilon}, \widetilde{u}^\varepsilon\right), \\ \widetilde{u}^\varepsilon(0,x) = \varepsilon^\gamma A_-(\psi_-)\left(\frac{-1}{\varepsilon}, \frac{x}{\varepsilon}\right). \end{cases}$$

If $\gamma > \alpha$, then the nonlinear term is negligible at leading order (as in Th. 2.7). If $\gamma = \alpha$, the nonlinear term is negligible away from $\{t=1\}$, and the caustic crossing is described in terms of the scattering operator $S$, as in Th. 2.5. This means that for a given equation, small initial data do not see the nonlinearity, and the critical size for the initial data corresponds to the case $\gamma = \alpha$. When $\gamma < \alpha$, strong nonlinear effects are expected, and we cannot give in general a description of the asymptotic behavior of $\widetilde{u}^\varepsilon$. For instance, it is proved in [9] that for a dissipative wave equation, the assumption $\gamma < \alpha$ leads to the absorption of the wave at the caustic (the wave is asymptotic to zero past the focus), while such a phenomenon is obviously impossible in the case of a conservative equation.

*Remark* 2.9. Theorem 2.7 implies the following asymptotics:

For every $0 < \tau < 1$, $\displaystyle\sup_{t \in [0,\tau]}\left\|u^\varepsilon(t,\cdot) - \frac{1}{(1-t)^\alpha}A_-(\psi_-)\left(\frac{t-1}{\varepsilon}, \frac{\cdot}{\varepsilon}\right)\right\|_{X,\varepsilon} \xrightarrow[\varepsilon \to 0]{} 0$,

For every $\tau > 1$, $\displaystyle\sup_{t \geq \tau}\left\|u^\varepsilon(t,\cdot) - \frac{1}{(t-1)^\alpha}A_+(\psi_-)\left(\frac{t-1}{\varepsilon}, \frac{\cdot}{\varepsilon}\right)\right\|_{X,\varepsilon} \xrightarrow[\varepsilon \to 0]{} 0$.

In this case, $u^\varepsilon$ behaves asymptotically like the free solution $w^\varepsilon$, for which it is known (see e.g. [12]) that the description of the focus crossing involves the Maslov index. Apparently, this index is missing in our writing: it is actually hidden in the definition of the operators $A_-$ and $A_+$. In the case of the Schrödinger equation for instance, we have

$$A_-(\varphi)(t,x) = e^{in\pi/4} e^{i\frac{x^2}{2t}} \widehat{\varphi}\left(\frac{x}{t}\right) \ , \quad A_+(\varphi)(t,x) = e^{-in\pi/4} e^{i\frac{x^2}{2t}} \widehat{\varphi}\left(\frac{x}{t}\right).$$

The presence of the Maslov index is measured by the phase shift $-n\pi/2$. This also means that in Theorem 2.5, two phenomena can be observed: the linear phenomenon, measured by the Maslov index (this is the transition $A_-/A_+$), and a nonlinear phenomenon, changing the profile of the solution ($\psi_-$ is replaced by $\psi_+ = S\psi_-$ past the focal point).

2.5. **A weakened version of Theorem 2.5.** It may happen that only the first two points of Assumption 1 are satisfied, i.e. that the group $U_0$ does not exhibit an asymptotic behavior in $\pm\infty$. This is for instance the case for the wave equation (see Section 5). In such a case, a weakened version of Theorem 2.5 still holds. Instead of the Cauchy problem (2.5), consider

$$(2.18) \quad \begin{cases} L(\varepsilon\partial)u^\varepsilon(t,x) = \varepsilon^{\alpha(p-1)} F\left(\dfrac{x}{\varepsilon}, u^\varepsilon(t,x)\right), \\ u^\varepsilon(0,x) = \varepsilon^{-\alpha} \left[U_0\left(\dfrac{-1}{\varepsilon}\right)\psi_-\right]\left(\dfrac{x}{\varepsilon}\right) + r^\varepsilon(x), \end{cases}$$

where $r^\varepsilon$ still denotes a $o(1)$ perturbation for the norm $\|\cdot\|_{X,\varepsilon}$. The only difference with (2.5) is that in the initial condition, $A_1(\psi_-)\left(\dfrac{-1}{\varepsilon}, \dfrac{x}{\varepsilon}\right)$ has been replaced by $\varepsilon^{-\alpha}\left[U_0\left(\dfrac{-1}{\varepsilon}\right)\psi_-\right]\left(\dfrac{x}{\varepsilon}\right)$. When the third point of Assumption 1 holds, the difference between these two expressions tends to 0 for the norm $\|\cdot\|_{X,\varepsilon}$ as $\varepsilon$ goes to 0. This fact is used in the first step of the proof of Th. 2.5 and is obviously superfluous if one considers (2.18) instead of (2.5). Similarly, if one is just interested in describing the caustic crossing in terms of the scattering operator $S$ without giving explicit asymptotics, then the third point of Assumption 1 is not needed in Steps 3 and 4 of the proof, and one has the following theorem:

**Theorem 2.10.** *Let $\psi_-, r^\varepsilon \in X$ and suppose that the first two points of Assumption 1 hold. Assume also that Assumptions 2 and 3 are satisfied. If $r^\varepsilon$ is asymptotically small, $\|r^\varepsilon\|_{X,\varepsilon} \xrightarrow[\varepsilon\to 0]{} 0$, then the following holds*

- *There exists $\varepsilon_0 > 0$ such that for all $0 < \varepsilon \leq \varepsilon_0$, there exists a unique solution $u^\varepsilon \in C(\mathbb{R}_+, X)$ to the Cauchy problem (2.18).*
- *For every $0 < \tau < 1$, $u^\varepsilon$ satisfies the following asymptotics,*

$$(2.19) \quad \sup_{t\in[0,\tau]} \left\| u^\varepsilon(t,\cdot) - \varepsilon^{-\alpha}\left[U_0\left(\dfrac{t-1}{\varepsilon}\right)\psi_-\right]\left(\dfrac{\cdot}{\varepsilon}\right) \right\|_{X,\varepsilon} \xrightarrow[\varepsilon\to 0]{} 0.$$

- *For every $\tau > 1$, $u^\varepsilon$ satisfies the following asymptotics,*

$$(2.20) \quad \sup_{t\geq\tau} \left\| u^\varepsilon(t,\cdot) - \varepsilon^{-\alpha}\left[U_0\left(\dfrac{t-1}{\varepsilon}\right)\psi_+\right]\left(\dfrac{\cdot}{\varepsilon}\right) \right\|_{X,\varepsilon} \xrightarrow[\varepsilon\to 0]{} 0,$$

*where $\psi_+ = S\psi_-$.*
- *There exists a "caustic profile". Let $\varphi = W_-\psi_- \in X$. Then*

$$\left\| u^\varepsilon(1,\cdot) - \dfrac{1}{\varepsilon^\alpha}\varphi\left(\dfrac{\cdot}{\varepsilon}\right) \right\|_{X,\varepsilon} \xrightarrow[\varepsilon\to 0]{} 0.$$

Similarly, Theorem 2.7 has an analog for problem (2.18), when the factor $\varepsilon^{\alpha(p-1)}$ is replaced by $\varepsilon^\delta$, with $\delta > \alpha(p-1)$. Turning the initial datum of $w^\varepsilon$ in (2.17) into that of (2.18), the statement is the same as in Theorem 2.7 and the proof is immediate.

2.6. **Relevance of Assumption 4.** Our point in this section is to show that in general, Assumption 4 is a consequence of the proof that shows that Assumptions 2 and 3 are satisfied. Suppose that Assumptions 2 and 3 are proved by a fixed point argument. It is the case in all the examples we consider. The Cauchy problem (2.4) with datum at $-\infty$,

$$L(\partial)\psi = F(x,\psi) \, , \quad U_0(-t)\psi(t)_{|t=-\infty} = \psi_-,$$

is solved by a fixed point argument applied to its integral formulation,

$$(2.21) \qquad \psi(t) = U_0(t)\psi_- + \int_{-\infty}^{t} U_0(t-s) F(x,\psi(s,x)) \, ds.$$

Assumption 2 means that the above integral is convergent, when $\psi_- \in X$, and Assumption 3 implies that this integral is defined for all $t \in \mathbb{R}$ and convergent as $t \to +\infty$.

To consider the small data case, fix $\psi_- \in X$, with $1 \le \|\psi_-\|_X \le 2$, and consider the problem with datum $\delta\psi_-$ for $\delta \in ]0,1]$,

$$L(\partial)\psi = F(x,\psi) \, , \quad U_0(-t)\psi(t)_{|t=-\infty} = \delta\psi_-.$$

Consider the new unknown $\widetilde{\psi} = \delta^{-1}\psi$. It solves

$$L(\partial)(\delta\widetilde{\psi}) = F(x,\delta\widetilde{\psi}) \, , \quad U_0(-t)\delta\widetilde{\psi}(t)_{|t=-\infty} = \delta\psi_-,$$

whose Duhamel's principle reads

$$\widetilde{\psi}(t) = U_0(t)\psi_- + \delta^{-1} \int_{-\infty}^{t} U_0(t-s) F\left(x, \delta\widetilde{\psi}(s,x)\right) ds.$$

Since we assume that $F(x,\cdot)$ is $p$-homogeneous for any $x \in \mathbb{R}^n$, we also have

$$\widetilde{\psi}(t) = U_0(t)\psi_- + \delta^{p-1} \int_{-\infty}^{t} U_0(t-s) F\left(x, \widetilde{\psi}(s,x)\right) ds.$$

The proof that Assumptions 2 and 3 are satisfied implies that there exists a constant $C$ independent of $\delta \in ]0,1]$ such that for any $t \in \mathbb{R}$,

$$\left\|\widetilde{\psi}(t) - U_0(t)\psi_-\right\|_X \le C\delta^{p-1}.$$

Back to $\psi$, we have

$$\|\psi(t) - U_0(t)\delta\psi_-\|_X \le C\delta^p \le C\|\delta\psi_-\|_X^p,$$

which is Assumption 4 with $C(\lambda) \propto \lambda^{p-1}$.

2.7. **The case of long range scattering.** Assumptions 2 and 3 mean that a complete scattering theory is available. More precisely, they indicate that the nonlinearity $F$ is negligible for large times, which is a feature of *short range* scattering. Only partial results are available for nonlinear long range scattering. In the case of the one-dimensional cubic Schrödinger equation (Eq. (1.2) with $n = \sigma = 1$), it is known that the nonlinear dynamics and the free dynamics cannot be compared ([37], [38], [2]); if $\psi_- \in L^2(\mathbb{R})$ and $U_0(-t)\psi(t) - \psi_- \xrightarrow[t \to -\infty]{} 0$ in $L^2$, where $\psi$ solves

$$(2.22) \qquad i\partial_t\psi + \frac{1}{2}\partial_x^2\psi = |\psi|^2\psi,$$

then $\psi = \psi_- = 0$. Ozawa ([31], see also [4]) proved that given $\psi_-$ sufficiently smooth and small, one can find $\psi$ solution of (2.22) such that

$$\left\| \psi(t) - e^{iS^-(t)} U_0(t)\psi_- \right\|_{L^2} \underset{t \to -\infty}{\longrightarrow} 0 \,, \tag{2.23}$$

where the phase shift $S^-$ is defined by

$$S^-(t) = \left| \widehat{\psi_-}\left(\frac{x}{t}\right) \right|^2 \log|t|.$$

Conversely, Hayashi and Naumkin [18] have proved that if $\psi$ solves (2.22), with an initial datum sufficiently small and smooth, then there exists $\psi_-$ such that (2.23) holds.

If we want to use these results in the context of geometrical optics, we introduce

$$v^\varepsilon(t,x) = \frac{1}{\sqrt{\varepsilon}} \psi\left(\frac{t-1}{\varepsilon}, \frac{x}{\varepsilon}\right).$$

The function $\psi_{|t=0}$ is obviously a concentrating profile for $v^\varepsilon$ at time $t = 1$, and the result of Hayashi and Naumkin shows that if $\psi_{|t=0} = \varphi$ is sufficiently small and smooth, then in $L^2$,

$$v^\varepsilon(0,x) \underset{\varepsilon \to 0}{\sim} \widehat{\varphi}(-x) e^{-i\frac{x^2}{2\varepsilon}} e^{-i\left|\widehat{\psi_-}(-x)\right|^2 \log \varepsilon}.$$

This suggests to replace (1.1) by the Cauchy problem

$$\begin{cases} i\varepsilon \partial_t u^\varepsilon + \frac{1}{2}\varepsilon^2 \partial_x^2 u^\varepsilon = \varepsilon |u^\varepsilon|^2 u^\varepsilon, \\ u^\varepsilon_{|t=0} = e^{-i\frac{|x|^2}{2\varepsilon}} e^{-i|f(x)|^2 \log \varepsilon} f(x). \end{cases}$$

This is exactly what was done in [4], and the above argument shows that the introduction of the $\log \varepsilon$ factor in the phase is necessary to have a concentrating profile at the focal point. A similar introduction was also performed in [7] in the case of a wave equation; again, this is due to long range scattering.

## 3. Hartree equations

As a first application of our general approach, we consider the Hartree equation. In this case, the nonlinear equation (2.3) writes

$$\begin{cases} i\partial_t \psi + \frac{1}{2}\Delta \psi = \lambda \left( |x|^{-\gamma} * |\psi|^2 \right) \psi, \ (t,x) \in \mathbb{R}_+ \times \mathbb{R}^n, \\ \psi_{|t=0} = \varphi \,, \end{cases} \tag{3.1}$$

with $n \geq 2$, $0 < \gamma < n$ and $\lambda \in \mathbb{R}$. The linear part of (3.1) is the Schrödinger equation, and we already checked that Assumption 1 is satisfied with $X = L^2(\mathbb{R}^n)$, $\alpha = n/2$ and

$$A_\pm(\varphi)(t,x) = e^{\mp in\frac{\pi}{4}} e^{i\frac{x^2}{2t}} \widehat{\varphi}\left(\frac{x}{t}\right).$$

It is easy to prove that the asymptotics (2.2) also holds in

$$\Sigma := H^1 \cap \mathcal{F}(H^1) = \left\{ \varphi \in H^1(\mathbb{R}^n) \,;\, |x|\varphi \in L^2(\mathbb{R}^n) \right\},$$

provided that $\varphi \in \Sigma$. On the other hand, the linear group $U_0(t) = e^{i\frac{t}{2}\Delta}$ is not unitary on $\Sigma$ (see Remark 2.3). We cannot apply directly the results of Section 2, but it is easy to adapt them in a similar fashion as in [3]. The main point consists in introducing the following time-dependent norm,

$$\left\|f\right\|_{\Sigma,t,\varepsilon} := \|f\|_{L^2} + \|\varepsilon\nabla f\|_{H^1} + \left\|\left(\frac{x}{\varepsilon} + i(t-1)\nabla\right) f\right\|_{L^2}. \tag{3.2}$$

The last operator, $x/\varepsilon + i(t-1)\nabla_x$, is the Galilean operator $J(t) = x + it\nabla_x$, which satisfies $J(t) = U_0(t)xU_0(-t)$, adapted to the scaling of our framework.

It is proved in [19] (see also [11]) that Assumptions 2, 3 and 4 hold in $\Sigma$ provided that

- $4/3 < \gamma < \min(4, n)$ and $\lambda > 0$, or,
- $1 < \gamma < \min(4, n)$, $\lambda > 0$, and $X$ is a small ball centered at the origin in $\Sigma$.

Since the linear part of (3.1) is the Schrödinger equation, the initial data we are interested in are the same as in [3], namely,

$$u^\varepsilon_{|t=0} = f(x)e^{-i\frac{|x|^2}{2\varepsilon}}, \text{ with } f \in \Sigma.$$

We can now state our result.

**Theorem 3.1.** *Let $f \in \Sigma$, $\lambda > 0$ and $1 < \gamma < \min(4, n)$, and consider the initial value problem*

$$i\varepsilon\partial_t u^\varepsilon + \frac{1}{2}\varepsilon^2\Delta u^\varepsilon = \varepsilon^\alpha\left(|x|^{-\gamma} * |u^\varepsilon|^2\right)u^\varepsilon; \quad u^\varepsilon_{|t=0} = f(x)e^{-i\frac{|x|^2}{2\varepsilon}}. \tag{3.3}$$

• *If $\alpha > \gamma$, then for $0 < \varepsilon \ll 1$, (3.3) has a unique solution $u^\varepsilon \in C(\mathbb{R}_+; \Sigma)$. In addition,*

$$\sup_{t \geq 0}\left\|u^\varepsilon(t) - w^\varepsilon(t)\right\|_{\Sigma,t,\varepsilon} \xrightarrow[\varepsilon \to 0]{} 0,$$

*where $w^\varepsilon$ is the solution to the initial value problem*

$$i\varepsilon\partial_t w^\varepsilon + \frac{1}{2}\varepsilon^2\Delta w^\varepsilon = 0; \quad w^\varepsilon_{|t=0} = f(x)e^{-i\frac{|x|^2}{2\varepsilon}}.$$

• *If $\alpha = \gamma$, suppose in addition that $\gamma > 4/3$ or that $\|f\|_\Sigma \ll 1$. Then for $0 < \varepsilon \ll 1$, (3.3) has a unique solution $u^\varepsilon \in C(\mathbb{R}_+; \Sigma)$. It satisfies the following asymptotics:*

- *If $0 < \tau < 1$, then*

$$\sup_{0 \leq t \leq \tau}\left\|u^\varepsilon(t,x) - \frac{1}{(1-t)^{n/2}}f\left(\frac{x}{1-t}\right)e^{i\frac{|x|^2}{2\varepsilon(t-1)}}\right\|_{\Sigma,t,\varepsilon} \xrightarrow[\varepsilon \to 0]{} 0.$$

- *If $\tau > 1$, then*

$$\sup_{t \geq \tau}\left\|u^\varepsilon(t,x) - \frac{e^{-in\pi/2}}{(t-1)^{n/2}}Zf\left(\frac{x}{1-t}\right)e^{i\frac{|x|^2}{2\varepsilon(t-1)}}\right\|_{\Sigma,t,\varepsilon} \xrightarrow[\varepsilon \to 0]{} 0,$$

*where*

$$Z = \mathcal{F}^{-1} \circ S \circ \mathcal{F},$$

*and $S$ is the nonlinear scattering operator associated to (3.1).*

## 4. KLEIN GORDON EQUATION

As a second application, we consider a nonlinear Klein-Gordon equation. Let $n \geq 2$. With the unknown $\psi = (\psi_1, \psi_2) : \mathbb{R}_+ \times \mathbb{R}^n \to \mathbb{R}^2$, the equation reads,

$$(4.1) \quad \begin{cases} \partial_t \psi + \begin{pmatrix} 0 & -1 \\ 1 - \Delta & 0 \end{pmatrix} \psi + \begin{pmatrix} 0 \\ \lambda |\psi_1|^{p-1} \psi_1 \end{pmatrix} = 0, \quad (t,x) \in \mathbb{R}_+ \times \mathbb{R}^n, \\ \psi|_{t=0} = \varphi, \end{cases}$$

with $\lambda \geq 0$, $p > 3$ if $n = 2$ and $1 + 4/n < p \leq (n+2)/(n-2)$ if $n \geq 3$.
This equation is obviously of the form (2.3) and is the usual defocusing nonlinear Klein-Gordon equation with power nonlinearity when written in scalar form,

$$\begin{cases} \partial_t^2 \psi_1 - \Delta \psi_1 + \psi_1 + \lambda |\psi_1|^{p-1} \psi_1 = 0, \\ \psi_1|_{t=0} = \varphi_1, \quad \partial_t \psi_1|_{t=0} = \varphi_2, \end{cases}$$

where $\psi = (\psi_1, \psi_2)$ and $\varphi = (\varphi_1, \varphi_2)$.

The natural space to work with is the energy space $X = H^1(\mathbb{R}^n) \times L^2(\mathbb{R}^n)$. From [15], [25], [28] and [29], we know that Assumptions 2 and 3 are satisfied. Concerning the asymptotic behavior of the unitary group $U_0$ associated to the linear part of (4.1), we have the following result, due to Nelson [30],

$$(4.2) \quad \left\| U_0(t)\varphi - \frac{1}{|t|^{n/2}} A_{\pm}(\varphi)(t, \cdot) \right\|_{L^2} \underset{t \to \pm\infty}{\longrightarrow} 0, \quad \text{for all } \varphi \in X,$$

where $A_{\pm} = (A_{1,\pm}, A_{2,\pm})$ is defined as $A_{\pm}(\varphi)(t,x) = 0$ for $|x| > |t|$ and, for $|x| < |t|$ as

$$(4.3) \quad \begin{aligned} A_{j,+}(\varphi)(x,t) &= \rho(t,x) \sum_{\pm} e^{\pm i\theta(t,x)} a_j^{\pm}\left(\frac{\mp x}{\sqrt{t^2 - |x|^2}}\right), \\ A_{j,-}(\varphi)(x,t) &= \rho(t,x) \sum_{\pm} e^{\mp i\theta(t,x)} a_j^{\pm}\left(\frac{\pm x}{\sqrt{t^2 - |x|^2}}\right), \end{aligned}$$

where $j = 1, 2$ and

$$(4.4) \quad \begin{aligned} \theta(t,x) &:= n\frac{\pi}{4} + \sqrt{t^2 - |x|^2}, \quad \rho(t,x) = \frac{1}{2}|t|(t^2 - |x|^2)^{-\frac{n+2}{4}}, \\ a_1^{\pm}(\xi) &= \widehat{\varphi_1}(\xi) \mp i \frac{1}{\sqrt{1 + |\xi|^2}} \widehat{\varphi_2}(\xi), \quad a_2^{\pm}(\xi) = \widehat{\varphi_2}(\xi) \pm i\sqrt{1 + |\xi|^2} \widehat{\varphi_1}(\xi). \end{aligned}$$

*Remark* 4.1. With the above definition of $A_{\pm}$, the initial conditions $A_{\pm}(\psi_-)(\frac{-1}{\varepsilon}, \frac{\cdot}{\varepsilon})$ for the Cauchy problem (2.5) are of size $O(1)$. For instance one has, for $\psi_-$ in the Schwartz class and $|x| < 1$,

$$A_{1,-}(\psi_-)\left(\frac{-1}{\varepsilon}, \frac{x}{\varepsilon}\right) = \frac{1}{2(1-|x|^2)^{\frac{n+2}{4}}} \sum_{\pm} e^{\mp i \frac{n\pi}{4}} e^{\mp i \frac{\sqrt{1-|x|^2}}{\varepsilon}} a_1^{\pm}\left(\frac{\mp x}{\sqrt{1-|x|^2}}\right)$$

Notice that the mapping $(\varphi_1, \varphi_2) \mapsto (a_1^-, a_1^+)$ is bijective on the Schwartz space $\mathcal{S}(\mathbb{R}^n)^2$, hence $a_1^{\pm}$ can be *any* functions in $\mathcal{S}(\mathbb{R}^n)$. On the other hand, $a_2^{\pm}$ and $a_1^{\pm}$ have to be related by (4.4).

The estimate (4.2) does not hold in general in the energy space $X$, and Assumption 1 is no more valid. However, we can prove:

**Lemma 4.2.** *Suppose $\varphi = (\varphi_1, \varphi_2)$ with $\varphi_1, \varphi_2$ in the Schwartz class $\mathcal{S}(\mathbb{R}^n)$. Then*

$$\left\| U_0(t)\varphi - \frac{1}{|t|^{n/2}} A_\pm(\varphi)(t,\cdot) \right\|_X \xrightarrow[t \to \pm\infty]{} 0.$$

*Proof.* The same methods that yield (4.2) also give

$$\left\| \nabla_x \left[U_0(t)\varphi\right] - \frac{1}{|t|^{n/2}} B_\pm(\varphi)(t,\cdot) \right\|_{L^2} \xrightarrow[t \to \pm\infty]{} 0,$$

where $B_\pm$ is defined as $A_\pm$ except that one has to replace $a_j^\pm(\xi)$ by $i\xi a_j^\pm(\xi)$, for $j = 1, 2$.

The lemma is therefore proved if we have

$$\frac{1}{|t|^{n/2}} \left\| \nabla_x \left[A_\pm(\varphi)(t)\right] - B_\pm(\varphi)(t) \right\|_{L^2} \xrightarrow[t \to \pm\infty]{} 0.$$

Under the assumption that $\varphi_1, \varphi_2$ are in the Schwartz class $\mathcal{S}(\mathbb{R}^n)$, easy computations show that it is the case. □

It follows from Lemma 4.2 that a weakened version of Assumption 1 holds. More precisely, the assumption remains true for all $\varphi$ in the Schwartz class instead of the energy space $X$, with $\alpha = n/2$ and $A_\pm$ as defined in (4.3)-(4.4).

It is easy to check that if the vectors $\psi_-$ and $r^\varepsilon$ are in $\mathcal{S}(\mathbb{R}^n)^2$, this weakened assumption permits the proof of the first two points of Th. 2.5, as well as of the last one. For the third point of the theorem, this is no longer true, because we do not know whether $\psi_+ := S\psi_-$ belongs to $\mathcal{S}(\mathbb{R}^n)^2$. We cannot use Lemma 4.2 in this case, but still can use the $L^2$ asymptotics of (4.2). This is the reason why, in the following theorem, the estimates after the caustic are given in norm $\|\cdot\|_{L^2} = \|\cdot\|_{L^2,\varepsilon}$ instead of $\|\cdot\|_{X,\varepsilon}$. Similarly, Assumption 4 is satisfied if $\psi_- \in \mathcal{S}(\mathbb{R}^n)^2$, and hence a weakened version of Th. 2.7 holds.

**Theorem 4.3.** *Let $\psi_- \in \mathcal{S}(\mathbb{R}^n)^2$, $\lambda > 0$, $n \geq 2$ and $1 + 4/n < p < 2^* - 1$, and consider the Cauchy problem*

(4.5)
$$\begin{cases} \varepsilon^2 \left(\partial_t^2 - \Delta\right) u^\varepsilon + u^\varepsilon + \lambda \varepsilon^{\alpha(p-1)} |u^\varepsilon|^{p-1} u^\varepsilon = 0, \\ (u^\varepsilon, \varepsilon \partial_t u^\varepsilon)|_{t=0} = A_-(\psi_-)\left(\frac{-1}{\varepsilon}, \frac{\cdot}{\varepsilon}\right), \end{cases}$$

*where $A_-(\psi_-)$ is given by (4.3)-(4.4).*

• *If $\alpha > n/2$ then for $0 < \varepsilon \ll 1$, (4.5) has a unique solution $u^\varepsilon \in C(\mathbb{R}_+; X)$, with $X = H^1 \times L^2$. In addition,*

$$\sup_{t \geq 0} \left\| u^\varepsilon(t) - w^\varepsilon(t) \right\|_{X,\varepsilon} \xrightarrow[\varepsilon \to 0]{} 0,$$

*where $w^\varepsilon$ is the solution to the initial value problem*

$$\begin{cases} \varepsilon^2 \left(\partial_t^2 - \Delta\right) w^\varepsilon + w^\varepsilon + \lambda \varepsilon^{\alpha(p-1)} |w^\varepsilon|^{p-1} w^\varepsilon = 0, \\ (w^\varepsilon, \varepsilon \partial_t w^\varepsilon)|_{t=0} = A_-(\psi_-)\left(\frac{-1}{\varepsilon}, \frac{\cdot}{\varepsilon}\right). \end{cases}$$

• *If $\alpha = n/2$, then for $0 < \varepsilon \ll 1$, (4.5) has a unique solution $u^\varepsilon \in C(\mathbb{R}_+; X)$. It satisfies the following asymptotics, where $A_\pm$ are defined in (4.3)-(4.4):*

- For every $0 < \tau < 1$,

(4.6)
$$\sup_{t\in[0,\tau]} \left\| u^\varepsilon(t,\cdot) - \frac{1}{(1-t)^{n/2}} A_{1,-}(\psi_-)\left(\frac{t-1}{\varepsilon}, \frac{\cdot}{\varepsilon}\right) \right\|_{L^2} +$$
$$\sup_{t\in[0,\tau]} \left\| \varepsilon\nabla_x u^\varepsilon(t,\cdot) - \frac{1}{(1-t)^{n/2}} \nabla_x A_{1,-}(\psi_-)\left(\frac{t-1}{\varepsilon}, \frac{\cdot}{\varepsilon}\right) \right\|_{L^2} +$$
$$\sup_{t\in[0,\tau]} \left\| \varepsilon\partial_t u^\varepsilon(t,\cdot) - \frac{1}{(1-t)^{n/2}} A_{2,-}(\psi_-)\left(\frac{t-1}{\varepsilon}, \frac{\cdot}{\varepsilon}\right) \right\|_{L^2} \xrightarrow[\varepsilon \to 0]{} 0 .$$

- For every $\tau > 1$,

(4.7)
$$\sup_{t\geq\tau} \left\| u^\varepsilon(t,\cdot) - \frac{1}{(t-1)^{n/2}} A_{1,+}(\psi_+)\left(\frac{t-1}{\varepsilon}, \frac{\cdot}{\varepsilon}\right) \right\|_{L^2} +$$
$$\sup_{t\geq\tau} \left\| \varepsilon\partial_t u^\varepsilon(t,\cdot) - \frac{1}{(t-1)^{n/2}} A_{2,+}(\psi_+)\left(\frac{t-1}{\varepsilon}, \frac{\cdot}{\varepsilon}\right) \right\|_{L^2} \xrightarrow[\varepsilon \to 0]{} 0 ,$$

where $\psi_+ = S\psi_-$.

## 5. WAVE EQUATION

We consider in this section a nonlinear wave equation in space dimension $n = 3$. In matricial form, with the unknown $\psi = (\psi_1, \psi_2) : \mathbb{R}_+ \times \mathbb{R}^3 \to \mathbb{R}^2$, it reads,

(5.1)
$$\begin{cases} \partial_t \psi + \begin{pmatrix} 0 & -1 \\ -\Delta & 0 \end{pmatrix} \psi + \begin{pmatrix} 0 \\ \lambda|\psi_1|^{p-1}\psi_1 \end{pmatrix} = 0, & (t,x) \in \mathbb{R}_+ \times \mathbb{R}^3, \\ \psi|_{t=0} = \varphi , \end{cases}$$

with $\lambda \geq 0$, and $p^*(3) < p \leq 5$, where
$$p^*(3) := \frac{5+\sqrt{33}}{4} .$$

This equation is obviously of the form (2.3) and is the usual defocusing nonlinear wave equation with power nonlinearity when written in scalar form,

$$\begin{cases} \partial_t^2 \psi_1 - \Delta\psi_1 + \lambda|\psi_1|^{p-1}\psi_1 = 0 , \\ \psi_1|_{t=0} = \varphi_1, \quad \partial_t\psi_1|_{t=0} = \varphi_2 , \end{cases}$$

where $\varphi = (\varphi_1, \varphi_2)$.

The natural space to work with here is the energy space $X = \dot{H}^1(\mathbb{R}^3) \times L^2(\mathbb{R}^3)$. It is known that Assumptions 2 and 3 are satisfied under the following assumption.

**Assumption 5.** *We assume $\lambda > 0$ and $p^*(3) < p \leq 5$. Denote*
$$\Sigma := \left\{ (f,g) \in H^1(\mathbb{R}^3) \times L^2(\mathbb{R}^3) \; ; \; |x|\nabla f, |x|g \in L^2(\mathbb{R}^3) \right\} .$$

*For $\psi_- \in X = \dot{H}^1(\mathbb{R}^3) \times L^2(\mathbb{R}^3)$, we suppose that:*

- *Either $p^*(3) < p < 5$ and $\psi_- \in \Sigma$,*
- *or $p^*(3) < p < 5$ and $\|\psi_-\|_X \ll 1$,*
- *or $p = 5$ and $\psi_- \in X$ is arbitrary.*

Assumption 5 implies Assumptions 2 and 3. For the first point, this was proven in [14] and [20]. For the second, we refer to [16]. Finally, the critical case $p = 5$ (last point) was studied in [35], [36], [25] and [1].

It is classical that the first two points of Assumption 1 are also satisfied with $\alpha = (n-1)/2 = 1$. In order to state the result concerning the asymptotic behavior of the group $U_0$, introduce the Radon transform $R$ defined for all $f$ in the Schwartz class $\mathcal{S}(\mathbb{R}^3)$ as

$$(Rf)(s,\omega) = \int_{\langle x,\omega\rangle = s} f(x) dS_x,$$

where $dS_x$ is the Euclidean measure on the hyperplane $\langle x, \omega \rangle = s$.
Let us also introduce the map

(5.2) $$k : \begin{cases} \mathcal{S}(\mathbb{R}^3) \times \mathcal{S}(\mathbb{R}^3) \to L^2(\mathbb{R}_s \times S^2_\omega), \\ (\varphi_1, \varphi_2) \mapsto \dfrac{1}{4\pi}\left(-\partial_s^2 R\varphi_1 + \partial_s R\varphi_2\right). \end{cases}$$

It is known (see e.g. [26], [32]) that $\mathcal{R}$ can be extended to $X$ onto $L^2(\mathbb{R} \times S^2)$ and that

(5.3) $$\left\| [U_0(t)\varphi]_2 \pm \frac{1}{|t|} k(\varphi)\left(|x| \mp |t|, \frac{x}{|x|}\right) \right\|_{L^2(\mathbb{R}^3)} \xrightarrow[t \to \pm\infty]{} 0,$$

and

(5.4) $$\left\| \nabla_x [U_0(t)\varphi]_1 \mp \frac{1}{|t|} k(\varphi)\left(|x| \mp |t|, \frac{x}{|x|}\right) \frac{x}{|x|} \right\|_{L^2(\mathbb{R}^3)} \xrightarrow[t \to \pm\infty]{} 0,$$

where $U_0(t)\varphi := ([U_0(t)\varphi]_1, [U_0(t)\varphi]_2)$.
However, these two estimates are not sufficient to ensure that the third point of Assumption 1 is satisfied. Indeed, if (5.4) gives the asymptotics of $\nabla u(t)$ for $|t|$ large, it does not tell us whether this asymptotic equivalent is the gradient of an element of $\dot{H}^1(\mathbb{R}^3)$. This is the reason why we can only use the "weakened" version of our results, i.e. Th. 2.10. Note that (5.3)-(5.4) can nevertheless be used to give asymptotics before and after the caustic. More precisely, we have:

**Theorem 5.1.** *Let $\psi_- \in X$, and consider the Cauchy problem*

(5.5) $$\begin{cases} \varepsilon^2 \left(\partial_t^2 - \Delta\right) u^\varepsilon + \lambda \varepsilon^{\alpha(p-1)} |u^\varepsilon|^{p-1} u^\varepsilon = 0, \\ (u^\varepsilon, \varepsilon \partial_t u^\varepsilon)|_{t=0} = \varepsilon^{-1} \left[U_0\left(\dfrac{-1}{\varepsilon}\right)\psi_-\right]\left(\dfrac{x}{\varepsilon}\right), \end{cases}$$

*where $U_0$ is the group associated to the free wave equation. Under Assumption 5, we have:*
• *If $\alpha > 1$, then for $0 < \varepsilon \ll 1$, (5.5) has a unique solution $u^\varepsilon \in C(\mathbb{R}_+; X)$, with $X = \dot{H}^1 \times L^2$. In addition,*

$$\frac{1}{\sqrt{\varepsilon}} \sup_{t \geq 0} \left\| \varepsilon \nabla_x u^\varepsilon(t) - \varepsilon \nabla_x w^\varepsilon(t) \right\|_{L^2} + \frac{1}{\sqrt{\varepsilon}} \sup_{t \geq 0} \left\| \varepsilon \partial_t u^\varepsilon(t) - \varepsilon \partial_t w^\varepsilon(t) \right\|_{L^2} \xrightarrow[\varepsilon \to 0]{} 0,$$

*where $w^\varepsilon$ is the solution to the initial value problem*

$$\begin{cases} \varepsilon^2 \left(\partial_t^2 - \Delta\right) w^\varepsilon + \lambda \varepsilon^{\alpha(p-1)} |w^\varepsilon|^{p-1} w^\varepsilon = 0, \\ (w^\varepsilon, \varepsilon \partial_t w^\varepsilon)|_{t=0} = \varepsilon^{-1} \left[U_0\left(\dfrac{-1}{\varepsilon}\right)\psi_-\right]\left(\dfrac{x}{\varepsilon}\right). \end{cases}$$

• *If $\alpha = 1$, then for $0 < \varepsilon \ll 1$, (5.5) has a unique solution $u^\varepsilon \in C(\mathbb{R}_+; X)$. It satisfies the following asymptotics, where $k(\psi_\pm)$ is defined in (5.2):*

- For every $0 < \tau < 1$,

$$\frac{1}{\sqrt{\varepsilon}} \sup_{t \in [0,\tau]} \left\| \varepsilon \partial_t u^\varepsilon(t, \cdot) - \frac{1}{(1-t)} k(\psi_-) \left( \frac{|x|+1-t}{\varepsilon}, \frac{x}{|x|} \right) \right\|_{L^2} \xrightarrow[\varepsilon \to 0]{} 0,$$

$$\frac{1}{\sqrt{\varepsilon}} \sup_{t \in [0,\tau]} \left\| \varepsilon \nabla_x u^\varepsilon(t, \cdot) + \frac{1}{(1-t)} k(\psi_-) \left( \frac{|x|+1-t}{\varepsilon}, \frac{x}{|x|} \right) \frac{x}{|x|} \right\|_{L^2} \xrightarrow[\varepsilon \to 0]{} 0.$$

- For every $\tau > 1$,

$$\frac{1}{\sqrt{\varepsilon}} \sup_{t \geq \tau} \left\| \varepsilon \partial_t u^\varepsilon(t, \cdot) + \frac{1}{(t-1)} k(\psi_+) \left( \frac{|x|-t-1}{\varepsilon}, \frac{x}{|x|} \right) \right\|_{L^2} \xrightarrow[\varepsilon \to 0]{} 0,$$

$$\frac{1}{\sqrt{\varepsilon}} \sup_{t \geq \tau} \left\| \varepsilon \nabla_x u^\varepsilon(t, \cdot) - \frac{1}{(t-1)} k(\psi_+) \left( \frac{|x|-t-1}{\varepsilon}, \frac{x}{|x|} \right) \frac{x}{|x|} \right\|_{L^2} \xrightarrow[\varepsilon \to 0]{} 0,$$

where $\psi_+ = S\psi_-$.

- There exists a caustic profile $(\varphi_1, \varphi_2) \in \dot{H}^1 \times L^2$ such that

$$\frac{1}{\sqrt{\varepsilon}} \left\| \varepsilon \nabla_x u^\varepsilon(1, \cdot) - \frac{1}{\varepsilon} \nabla_x \varphi_1 \left( \frac{\cdot}{\varepsilon} \right) \right\|_{L^2} \xrightarrow[\varepsilon \to 0]{} 0,$$

$$\frac{1}{\sqrt{\varepsilon}} \left\| \varepsilon \partial_t u^\varepsilon(1, \cdot) - \frac{1}{\varepsilon} \varphi_2 \left( \frac{\cdot}{\varepsilon} \right) \right\|_{L^2} \xrightarrow[\varepsilon \to 0]{} 0.$$

*Remark* 5.2. In this theorem, we estimate the $\varepsilon$-derivatives of $u^\varepsilon$. Since we adopt the viewpoint of geometrical optics, this must not be too surprising, as explained in Section 2.3. We restate our result with a different viewpoint below.

In the above asymptotics, the factor $\varepsilon^{-1/2}$ may suggest that we have approximations with an error term of order $o(\sqrt{\varepsilon})$, which is much more precise than $o(1)$. We must bear in mind that each term whose $L^2$-norm is assessed is actually of order $O(\sqrt{\varepsilon})$ in $L^2$, so our result is to be compared to what is usually an $o(1)$-asymptotics.

As noticed in Remark 2.2, the initial data chosen in (5.5) are of order $O(1)$ in $L^\infty$. More usual in the study of nonlinear wave equations is the energy norm,

$$E^\varepsilon_\lambda(t) = \frac{1}{2} \int |\partial_t u^\varepsilon(t,x)|^2 dx + \frac{1}{2} \int |\nabla_x u^\varepsilon(t,x)|^2 dx + \frac{\lambda \varepsilon^{\alpha(p-1)-2}}{p+1} \int |u^\varepsilon(t,x)|^{p+1} dx.$$

If the initial datum is in the energy space $X$, $E^\varepsilon_\lambda$ does not depend on time when $u^\varepsilon$ solves (5.5). In particular, $E^\varepsilon_0$ is the energy of the free wave equation. When we consider initial data as in (5.5), then with $\alpha \geq 1$ and $p^*(3) < p \leq 5$, the initial free energy $E^\varepsilon_0$, as well as the nonlinear energy $E^\varepsilon_\lambda$, are of order $O(\varepsilon^{-1})$. Modulating the size of $u^\varepsilon$ by a factor $\varepsilon^\gamma$ like we did in Remark 2.8, we can alter the initial value problem (5.5) so that the energy is of order $O(1)$. Define $\mathfrak{u}^\varepsilon(t,x) := \sqrt{\varepsilon} u^\varepsilon(t,x)$. Then (5.5) is equivalent to

$$(5.6) \quad \begin{cases} \left( \partial_t^2 - \Delta \right) \mathfrak{u}^\varepsilon + \lambda \varepsilon^\beta |\mathfrak{u}^\varepsilon|^{p-1} \mathfrak{u}^\varepsilon = 0, \\ (\varepsilon^{-1} \mathfrak{u}^\varepsilon, \partial_t \mathfrak{u}^\varepsilon)|_{t=0} = \varepsilon^{-3/2} \left[ U_0 \left( \frac{-1}{\varepsilon} \right) \psi_- \right] \left( \frac{x}{\varepsilon} \right), \end{cases}$$

with $\beta = \alpha(p-1) - p/2 - 3/2$. In particular, in the case of the critical nonlinearity, $p = 5$, with critical scaling as far as focusing is concerned, $\alpha = 1$, we have $\beta = 0$. That is, $\mathfrak{u}^\varepsilon$ satisfies the critical wave equation with no $\varepsilon$,

$$(5.7) \quad \left( \partial_t^2 - \Delta \right) \mathfrak{u}^\varepsilon + \lambda |\mathfrak{u}^\varepsilon|^4 \mathfrak{u}^\varepsilon = 0,$$

and its energy is bounded, in $O(1)$. This model meets the more general problem studied in [1], which consists in characterizing the bounded energy sequences of solutions to the above equation. Notice that it that statement, no scale $\varepsilon$ is privileged: several scales may have a leading order influence. It is proved in that paper that nonlinear effects may occur, measured by a scattering operator as in Theorem 5.1. The authors prove even more general results, some of which we examine in Section 6.

Now consider the case of a subcritical nonlinearity, as far as the existence of solutions in the energy space is concerned, that is $p < 5$, and for $\gamma \in \mathbb{R}$, introduce the initial value problem

(5.8) $$\begin{cases} \left(\partial_t^2 - \Delta\right) u^\varepsilon + \lambda \varepsilon^\gamma |u^\varepsilon|^{p-1} u^\varepsilon = 0, \\ (\varepsilon^{-1} u^\varepsilon, \partial_t u^\varepsilon)|_{t=0} = \varepsilon^{-3/2} \left[U_0\left(\frac{-1}{\varepsilon}\right)\psi_-\right]\left(\frac{x}{\varepsilon}\right), \end{cases}$$

It is natural to introduce a critical value for $\gamma$, which corresponds to the critical case $\alpha = 1$ and $\gamma = \beta$,
$$\gamma_c = \frac{p-5}{2} < 0 \ .$$
If $\gamma > \gamma_c$, then solutions to (5.8) are bounded energy sequences, which are asymptotically linear. This is the case in particular if $\gamma = 0$, and we retrieve an aspect of the results stated in [13]. In this paper, the author prove that *any* bounded energy sequences of initial data lead to solutions to the subcritical wave equation which are asymptotically linear ("linearizable"). Those initial data are much more general than ours, but since ours are those that cause focusing at one point, and are therefore the "worst" initial data, it seems that our result are in a way, related.

On the other hand, if $\gamma = \gamma_c$, nonlinear effects are relevant at the focal point, and are measured by the scattering operator $S$. Since $\gamma_c < 0$, this means that an amplifying term is necessary for the nonlinearity to act on solutions whose energy is $O(1)$. Another viewpoint consists in considering solutions to an equations with no $\varepsilon$, whose energy depends on $\varepsilon$. Considering
$$\widetilde{u}^\varepsilon = \varepsilon^{\gamma_c/(p-1)} u^\varepsilon = \varepsilon^{\frac{p-5}{2(p-1)}} u^\varepsilon = \varepsilon^{\frac{p-3}{p-1}} u^\varepsilon,$$
we have $\widetilde{u}^\varepsilon$ solution of

(5.9) $$\left(\partial_t^2 - \Delta\right)\widetilde{u}^\varepsilon + \lambda |\widetilde{u}^\varepsilon|^{p-1}\widetilde{u}^\varepsilon = 0,$$

and $E_\lambda(\widetilde{u}^\varepsilon) = O(\varepsilon^{(p-5)/(p-1)})$. This is therefore the minimal size of the energy for $\varepsilon$-oscillatory solutions to (5.9) not to be linearizable.

*Remark* 5.3. The same discussion would give similar results for the Klein-Gordon equation (4.5), with the difference that $\varepsilon$ will always be present in the equation.

To conclude this section, we apply our results to a semilinear wave equation in space dimension three, that is not of the form (5.5), but

(5.10) $$\begin{cases} \Box u^\varepsilon + a\,\varepsilon^\gamma |\partial_t u^\varepsilon|^{p-1}\,\partial_t u^\varepsilon = 0, \\ u^\varepsilon\big|_{t=0} = \varepsilon U_0\left(r, \frac{r-r_0}{\varepsilon}\right), \qquad \partial_t u^\varepsilon\big|_{t=0} = U_1\left(r, \frac{r-r_0}{\varepsilon}\right), \end{cases}$$

where $r = |x|$ and $a \in \mathbb{C}$. This equation was studied in [8], [10] and [9], for initial data $U_j$ which are smooth, and compactly supported in their last variable; this means that the data are *pulse like* (as opposed to wave trains, for which the $U_j$'s are periodic in their last argument). Since the initial data are spherical, so is the solution. The changes of unknowns

$$(5.11) \qquad \tilde{u}^\varepsilon(t,r) := r u^\varepsilon(t,r), \qquad v^\varepsilon_\mp := (\partial_t \pm \partial_r)\tilde{u}^\varepsilon, \qquad v^\varepsilon_\mp \in C^\infty(\mathbb{R}_t \times \mathbb{R}_r),$$

turns the initial 3D problem (5.10) into a one-dimensional mixed problem,

$$(5.12) \quad \begin{cases} (\partial_t \pm \partial_r) v^\varepsilon_\pm = \varepsilon^\gamma r^{1-p} g(v^\varepsilon_- + v^\varepsilon_+), & r > 0, \qquad g(y) := -a 2^{-p} |y|^{p-1} y, \\ (v^\varepsilon_- + v^\varepsilon_+)\big|_{r=0} = 0, \\ v^\varepsilon_{\pm | t=0} = v^\varepsilon_{\pm 0}, \end{cases}$$

where $v^\varepsilon_{\pm 0}, \varepsilon \partial_t v^\varepsilon_{\pm 0} \in L^\infty(\mathbb{R}_+)$. The linear part of (5.12) is the free wave equation in space dimension one in the quarter of plane $t, r \geq 0$. It is easy to see that Assumption 1 is satisfied, with $X = L^\infty(\mathbb{R} \times \mathbb{R}_+)^2$ and $\alpha = 0$. The main result of [10] consists in proving Theorem 2.5 for (5.12), which corresponds to the case $\gamma = p - 2 > 0$. The core of the paper is actually the proof that Assumptions 2 and 3 are satisfied, and demands to construct the scattering operator. As mentioned in Section 2.6, deducing Theorem 2.7 is then easy, this is a part of [9]; if $\gamma > p - 2 > 0$, then the solutions to (5.12) are asymptotically free. They are also asymptotically free in the case $\alpha > 0$ and $1 < p \leq 2$, but the proof is slightly different, and the result is weaker in some sense: the time estimates are local, while they are global in the case $\gamma > p - 2 > 0$.

## 6. Remarks and comments

To end this paper, we outline the limits of our results. Nevertheless, we believe they are an interesting generalization, and provide a better understanding, of the results of [3] and [10].

6.1. **Geometry.** First, the geometrical background is very simple; we consider a caustic reduced to a point. This is at least a step in the understanding of nonlinear phenomena at a caustic in geometrical optics. More general geometries are considered in [24], but nonlinear phenomena at the caustic are well understood only in the case of dissipative equations, leading to absorption of oscillations. We emphasized critical indexes in the case of a linear geometric optics régime, with a focal point, and gave a description of the critical case. Understanding what may happen in the supercritical case ($\delta < \alpha(p-1)$ in (2.16)) should lead to interesting phenomena.

6.2. **Other focal points.** We consider initial data that focus at the origin at time $t = 1$ (see Section 2.3). It is easy to modify our approach so that focusing occurs at $(t, x) = (t^\varepsilon_j, x^\varepsilon_j)$ for any prescribed $(t^\varepsilon_j, x^\varepsilon_j)$. The question to know what happens when a superposition of such initial data (taking different pairs $(t^\varepsilon_j, x^\varepsilon_j)$) is natural. Answers are available in the cases of the critical wave equation (5.7) and of the semiclassical nonlinear Schrödinger equation (1.1), given in [1] and [5] (see also [6]) respectively. In these papers, nonlinear superposition principles are established. Essentially, the nonlinearity is negligible when no focusing occurs, for the régimes considered of those of *linear* geometric optics. When focusing occurs, it does at a scale which is so small ($\varepsilon$ in the case of (1.1)) that if the pairs $(t^\varepsilon_j, x^\varepsilon_j)$ are chosen "orthogonal" in a suitable sense, then no interaction occurs.

FOCUSING AT A POINT WITH CAUSTIC CROSSING     23### 6.3. The role of $\varepsilon$.
We consider only one small parameter $\varepsilon > 0$, as is usual in geometrical optics. However, other small parameters might play a similar role. As we already mentioned, no specific scale appears in the study of bounded energy sequences solutions to the critical wave equation (5.7). Indeed, H. Bahouri and P. Gérard show that several small parameters, of different orders, are necessary to describe such solutions. This is because of the scaling of the nonlinearity, which is critical for the wave equation.

### 6.4. More general initial data.
The role of such initial data as those we consider in (2.5) or (2.18) could be better understood. Our argument is that these data cause focusing at one point, which is the most degenerate caustic. Therefore, the amplitude growth is the most important predicted by geometrical optics. Intuitively, this suggests that such data are the only ones, which are of the same order in the norm $\|\cdot\|_{X,\varepsilon}$ (that is, of order $O(1)$), that make the nonlinear relevant in (2.5) or (2.18). This means that solutions to Equation (2.5) whose initial data are bounded in the norm $\|\cdot\|_{X,\varepsilon}$ are not linearizable if and only if the initial data can be written as a superposition of such terms as those we consider, plus a term which generates a linearizable solution. This vague statement is made rigorous in the cases of the critical wave equation (5.7) and of the semiclassical nonlinear Schrödinger equation (1.1), in [1] and [5] respectively. The proofs of these results rely on estimates which use finer properties of Equations (2.1) and (2.3) than those we state in Assumptions 1, 2 and 3. The generality of these results still has to be understood.

### 6.5. Stability properties.
In (2.5), we can perturb the initial datum by a term $r^\varepsilon$ without changing the asymptotics of $u^\varepsilon$ provided that $\|r^\varepsilon\|_{X,\varepsilon} \to 0$ as $\varepsilon \to 0$. This is a stability property for (2.5). As we saw in the proof of Theorem 2.5, this stability property stems from a continuity property of the wave operators for (2.3) (Assumption 2), and from the global well-posedness for (2.3) (Assumption 3).

**Acknowledgements.** The authors want to thank J. Rauch for his fruitful comments. This work was partially supported by the ACI Jeunes Chercheurs du Ministère de la Recherche "Équation des ondes : oscillations, dispersion et contrôle" and "Solutions oscillantes d'EDP", and by GDR 2103 EAPQ CNRS.

## References

1. H. Bahouri and P. Gérard, *High frequency approximation of solutions to critical nonlinear wave equations*, Amer. J. Math. **121** (1999), no. 1, 131–175. MR **2000i:**35123
2. J. E. Barab, *Nonexistence of asymptotically free solutions for nonlinear Schrödinger equation*, J. Math. Phys. **25** (1984), 3270–3273.
3. R. Carles, *Geometric optics with caustic crossing for some nonlinear Schrödinger equations*, Indiana Univ. Math. J. **49** (2000), no. 2, 475–551. MR **2001k:**35265
4. ______, *Geometric optics and long range scattering for one-dimensional nonlinear Schrödinger equations*, Comm. Math. Phys. **220** (2001), no. 1, 41–67.
5. R. Carles, C. Fermanian, and I. Gallagher, *On the role of quadratic oscillations in nonlinear Schrödinger equations*, Preprint, `arXiv:math.AP/0212171`, 2002.
6. ______, *Rôle des oscillations quadratiques dans des équations de Schrödinger non linéaires*, Séminaire sur les Équations aux Dérivées Partielles, 2002–2003 (Palaiseau), École Polytech., 2002, Exp. No. IX.
7. R. Carles and D. Lannes, *Focusing of a pulse with arbitrary phase shift for a nonlinear wave equation*, Bull. Soc. Math. France (2003), To appear.

8. R. Carles and J. Rauch, *Focusing of spherical nonlinear pulses in $\mathbb{R}^{1+3}$*, Proc. Amer. Math. Soc. **130** (2002), no. 3, 791–804. MR 1 866 035
9. ______, *Focusing of Spherical Nonlinear Pulses in $\mathbb{R}^{1+3}$ III. Sub and Supercritical cases*, Preprint, `arXiv:math.AP/0212288`, 2002.
10. ______, *Focusing of Spherical Nonlinear Pulses in $\mathbb{R}^{1+3}$ II. Nonlinear Caustic*, Rev. Mat. Iberoamericana **20** (2004), no. 2, to appear.
11. T. Cazenave, *An introduction to nonlinear Schrödinger equations*, Text. Met. Mat., vol. 26, Univ. Fed. Rio de Jan., 1993.
12. J. J. Duistermaat, *Oscillatory integrals, Lagrange immersions and unfolding of singularities*, Comm. Pure Appl. Math. **27** (1974), 207–281. MR 53 #9306
13. P. Gérard, *Oscillations and concentration effects in semilinear dispersive wave equations*, J. Funct. Anal. **141** (1996), no. 1, 60–98. MR **97k:**35171
14. J. Ginibre and G. Velo, *Conformal invariance and time decay for nonlinear wave equations. I, II*, Ann. Inst. H. Poincaré Phys. Théor. **47** (1987), no. 3, 221–261, 263–276. MR **89e:**35021
15. ______, *The global Cauchy problem for the nonlinear Klein-Gordon equation. II*, Ann. Inst. H. Poincaré Anal. Non Linéaire **6** (1989), no. 1, 15–35. MR **90b:**35153
16. ______, *Scattering theory in the energy space for a class of nonlinear wave equations*, Comm. Math. Phys. **123** (1989), no. 4, 535–573. MR **90i:**35172
17. O. Guès, *Développement asymptotique de solutions exactes de systèmes hyperboliques quasilinéaires*, Asymptotic Anal. **6** (1993), no. 3, 241–269.
18. N. Hayashi and P. Naumkin, *Asymptotics for large time of solutions to the nonlinear Schrödinger and Hartree equations*, American Journal of Mathematics **120** (1998), 369–389.
19. N. Hayashi and Y. Tsutsumi, *Scattering theory for Hartree type equations*, Ann. Inst. H. Poincaré Phys. Théor. **46** (1987), no. 2, 187–213. MR **89a:**81158
20. K. Hidano, *Scattering problem for the nonlinear wave equation in the finite energy and conformal charge space*, J. Funct. Anal. **187** (2001), no. 2, 274–307. MR **2002k:**35242
21. J. Hunter and J. Keller, *Caustics of nonlinear waves*, Wave motion **9** (1987), 429–443.
22. S. Ibrahim, *Geometric optics for nonlinear concentrating waves in focusing and non-focusing two geometries*, Preprint, 2002.
23. J.-L. Joly, G. Métivier, and J. Rauch, *Generic rigorous asymptotic expansions for weakly nonlinear multidimensional oscillatory waves*, Duke Math. J. **70** (1993), no. 2, 373–404.
24. J.-L. Joly, G. Métivier, and J. Rauch, *Caustics for dissipative semilinear oscillations*, Mem. Amer. Math. Soc. **144** (2000), no. 685, viii+72. MR **2000i:**35115
25. L. Kapitanski, *Global and unique weak solutions of nonlinear wave equations*, Math. Res. Lett. **1** (1994), no. 2, 211–223. MR **95f:**35158
26. P. D. Lax and R. S. Phillips, *Scattering theory*, second ed., Pure and Applied Mathematics, vol. 26, Academic Press Inc., Boston, MA, 1989, With appendices by Cathleen S. Morawetz and Georg Schmidt. MR **90k:**35005
27. D. Ludwig, *Uniform asymptotic expansions at a caustic*, Comm. Pure Appl. Math. **19** (1966), 215–250.
28. K. Nakanishi, *Scattering theory for the nonlinear Klein-Gordon equation with Sobolev critical power*, Internat. Math. Res. Notices (1999), no. 1, 31–60. MR **2000a:**35174
29. ______, *Remarks on the energy scattering for nonlinear Klein-Gordon and Schrödinger equations*, Tohoku Math. J. (2) **53** (2001), no. 2, 285–303. MR **2002e:**35220
30. S. Nelson, *$L^2$ asymptotes for the Klein-Gordon equation*, Proc. Amer. Math. Soc. **27** (1971), 110–116. MR 42 #6444
31. T. Ozawa, *Long range scattering for nonlinear Schrödinger equations in one space dimension*, Comm. Math. Phys. **139** (1991), 479–493.
32. V. Petkov, *Scattering theory for hyperbolic operators*, Studies in Mathematics and its Applications, vol. 21, North-Holland Publishing Co., Amsterdam, 1989. MR **91e:**35170
33. J. Rauch, *Partial differential equations*, Graduate Texts in Math., vol. 128, Springer-Verlag, New York, 1991.
34. J. Rauch and M. Keel, *Lectures on geometric optics*, Hyperbolic equations and frequency interactions (Park City, UT, 1995), Amer. Math. Soc., Providence, RI, 1999, pp. 383–466.
35. J. Shatah and M. Struwe, *Regularity results for nonlinear wave equations*, Ann. of Math. (2) **138** (1993), no. 3, 503–518. MR **95f:**35164
36. ______, *Well-posedness in the energy space for semilinear wave equations with critical growth*, Internat. Math. Res. Notices (1994), no. 7, 303–309. MR **95e:**35132

37. W. A. Strauss, *Nonlinear scattering theory*, Scattering theory in mathematical physics (J. Lavita and J. P. Marchand, eds.), Reidel, 1974.
38. ______, *Nonlinear scattering theory at low energy*, J. Funct. Anal. **41** (1981), 110–133.

MAB, UMR 5466 CNRS, 351 cours de la Libération, 33 405 Talence cedex, France
*E-mail address*: `carles@math.u-bordeaux.fr`
*URL*: `www.math.u-bordeaux.fr/~carles`

MAB, UMR 5466 CNRS, 351 cours de la Libération, 33 405 Talence cedex, France
*E-mail address*: `lannes@math.u-bordeaux.fr`
*URL*: `www.math.u-bordeaux.fr/~lannes`